\documentstyle[12pt]{article}
\begin{document}
\title{On the Stability of $J^*$-Homomorphisms \footnote{
The first author was supported by Korea Research Foundation Grant
KRF-2004-041-C00023.}}
\author{Choonkil Baak and  Mohammad Sal Moslehian}
\date{}
\maketitle{}
\begin{abstract}
The main purpose of this paper is to prove the generalized Hyers-Ulam-Rassias stability of $J^*$-homomorphisms between $J^*$-algebras.
\end{abstract}

\section {Introduction}

In 1940, S. M. Ulam posed the following problem concerning the stability of group homomorphisms [13]:

Given a group $G_1$, a metric group $(G_2, d)$ and a positive number $\epsilon$, does there exist a $\delta>0$
such that if a function $f : G_1 \to G_2$ satisfies the inequality $d(f(xy), f(x)f(y))<\delta$
for all $x, y \in G_1$ then there exists a homomorphism $T : G_1 \to G_2$ such that $d(f(x), T(x))<\epsilon$ for all $x\in G_1$. If this problem has a solution, we say that the homomorphisms
from $G_1$ to $G_2$ are stable or the functional equation $f(xy)=f(x)f(y)$ is stable.

In 1941, D. H. Hyers gave a partial solution of Ulam's problem in the context of Banach spaces
as the following [7]:

Suppose that $E_1, E_2$ are Banach spaces and $f : E_1 \to E_2$ satisfies the following
condition: there is $\epsilon>0$ such that $\|f(x+y)-f(x)-f(y)\|<\epsilon$ for all $x, y \in E_1$. Then there is an additive mapping $T : E_1 \to E_2$ such that $\|f(x)-T(x)\|<\epsilon$ for all $x\in E_1$.

Now, assume that $E \sb 1$ and $E\sb 2$ are real normed spaces with $E\sb 2$ complete, $f: E\sb 1\to E\sb 2$ is a mapping such that $f(tx)$ is continuous in $t\in R$ for each fixed $x\in E\sb 1$, and let there exist $\varepsilon\ge 0$ and $p\in [0, \infty)-\{1\}$ such that
$$\|f(x+y)-f(x)-f(y)\|\le \varepsilon(\|x\|\sp p+\|y\|\sp p),$$ for all $x, y \in E\sb 1$.
It was shown by Th. M. Rassias [11] for $p\in [0, 1)$ and Z. Gajda [2] for $p>1$ that there exists a unique linear map $T: E\sb 1 \to E\sb 2$ such that $$\|f(x)-T(x)\|\le\frac{2\epsilon}{|2^p-2|}\|x\|\sp p,$$ for all $x \in E\sb 1$. Using Hyers' method, indeed,
$T(x)$ is defined by $\lim\sb{n\to\infty}2\sp{-n}f(2\sp nx)$ if $0\leq p<1$, and
$\lim\sb{n\to\infty}2\sp nf(2\sp{-n}x)$ if $p>1$. This phenomenon is called the Hyers--Ulam--Rassias stability. It is shown that there is no analogue of Rassias result for $p=1$ (see [12]).

In 1992, a generalization of Rassias' theorem was obtained by G\u avruta as follows [3]:

Suppose $(G,+)$ is an abelian group, $E$ is a Banach space and the so-called admissible
 control function $\varphi:G\times G\to [0, \infty)$ satisfies $$\widetilde{\varphi}(x, y):
=\frac{1}{2}\displaystyle{\sum_{n=0}^\infty}2^{-n}\varphi(2^n x,2^n y)<\infty,$$ for all $x,y\in G$. If $f : G \to E$ is a mapping with $$\|f(x+y)-f(x)-f(y)\|\leq \varphi(x,y),$$ for all $x, y \in G$, then there exists a unique mapping $T : G \to E$ such that $T(x+y)=T(x)+T(y)$ and $\|f(x)-T(x)\|\leq\widetilde{\varphi}(x, x)$ for all $x,y\in G$.

During the last decades several stability problems of functional equations have been investigated in spirit of Hyers--Ulam--Rassias, see [1].

By a $J^*$-algebra we mean a closed subspace $A$ of a $C^*$-algebra such that
$xx^*x\in A$ whenever $x\in A$. Many familiar spaces are $J^*$-algebras [4]. For example, {\rm (i)} every Cartan factor of type $I$, i.e.,  the space of all bounded operators $B(H,K)$ between Hilbert spaces $H$ and $K$;
{\rm (ii)} every Cartan factor of type $IV$, i.e.,  a closed $*$-subspace $A$ of $B(H)$ in which
the square of each operator in $A$ is a scalar multiple of the identity operator on $H$; {\rm (iii)} every $JC^*$-algebra; {\rm (iv)} every ternary algebra of operators [5].

A $J^*$-homomorphism between $J^*$-algebras $A$ and $B$ is defined to be a linear mapping $L : A \to B$ such that $T(xx^*x)=Tx(Tx)^*Tx$ for all $x\in A$. In particular, every $*$-homomorphism between $C^*$-algebras is a $J^*$-homomorphism.

In [9], the first author establishes the stability of $*$-homomorphisms of a
$C^*$-algebra. Using Rassias' technique [11] and some similar methods as
in [10], we generalize some of the results for $J^*$-homomorphisms.

Throughout this paper, $A$ and $B$ denote $J^*$-algebras.

\section {Main results}

Our main purpose is to prove the generalized Hyers--Ulam--Rassias stability of $J^*$-homomorphisms.

{\bf Proposition 2.1.} {\it Let $r>1$, and let $T : A \to A$ be a mapping satisfying $T(rx)=rT(x)$ for all $x\in A$ and let there exist
a function $\varphi: A\times A\times A\to [0, \infty)$ such that
$$\displaystyle{\lim_{n\to\infty}}r^{-n}\varphi(r^n x, r^n y, r^n z)=0,$$
\begin{equation}
\|T(\lambda x+y+zz^*z)-\lambda T(x)-T(y)-T(z)T(z)^*T(z)\|\leq\varphi(x, y, z),
\end{equation}
for all $\lambda \in {\bf C}$ and all $x, y, z\in A$.
Then $T$ is a $J^*$-homomorphism.}

{\bf Proof.} $T(0)=0$ since $T(0)=rT(0)$. Put $x=y=0$ in {\rm (1)}. Then
\begin{eqnarray*}
\|T(zz^*z)-T(z)T(z)^*T(z)\|&=& \frac{1}{r^{3n}}\|T(r^nz(r^nz)^*r^nz)-T(r^nz)T(r^nz)^*T(r^nz)\| \\
& \leq & \frac{1}{r^{3n}}\varphi(r^nz, r^nz, r^nz)  \leq  \frac{1}{r^{n}}\varphi(r^nz, r^nz, r^nz),
\end{eqnarray*}
for all $z \in A$.
The right side tends to zero as $n\to\infty$.
So  $$T(zz^*z)=T(z)T(z)^*T(z),$$ for all $z \in A$.

Similarly,
 one can shows that $$T(\lambda x+y)=\lambda T(x)+T(y),$$
for all $x, y \in A$. $\Box$

{\bf Theorem 2.2.} {\it Suppose $h : A \to B$ is a mapping with $h(0)=0$ for which there exists a function
$\varphi: A\times A\times A\to [0, \infty)$ such that
$$\widetilde{\varphi}(x, y, z):=\frac{1}{2}\displaystyle{\sum_{n=0}^\infty}2^{-n}\varphi(2^nx, 2^ny, 2^nz)<\infty,$$
\begin{equation}
\|h(\mu x+\mu y+zz^*z)-\mu h(x)-\mu h(y)-h(z)h(z)^*h(z)\|\leq \varphi(x, y, z),
\end{equation}
for all $\mu\in S^1=\{\lambda\in {\bf C} : |\lambda|=1\}$ and all $x, y, z\in A$.
Then there exists a unique $J^*$-homomorphism $T : A \to B$ such that $$\|h(x)-T(x)\|\leq\widetilde{\varphi}(x, x, 0),$$
for all $x\in A$.}

{\bf Proof.} In {\rm (2)}, assume that $z=0$ and $\mu=1$. Then the G\u avruta theorem implies that
there is a unique additive mapping $T : A \to B$ given by $$T(x)=\displaystyle{\lim_{n\to\infty}}\frac{h(2^nx)}{2^n},$$
for all $ x\in A$.
By {\rm (2)}, $$\|h(2^n\mu x)-2\mu h(2^{n-1}x)\|\leq\varphi(2^{n-1}x, 2^{n-1}x, 0),$$
for all $ x\in A$ and all $ \mu\in S^1$.
Then
$$\|\mu h(2^nx)-2\mu h(2^{n-1}x)\|\leq|\mu|\cdot \|h(2^nx)-2h(2^{n-1}x)\|\leq\varphi(2^{n-1}x, 2^{n-1}x, 0),$$
for all $ \mu\in S^1$ and all $x\in A$. So
\begin{eqnarray*}
\|2^{-n}h(2^n\mu x)-\mu 2^{-n}h(2^n x)\| & \leq & 2^{-n}\|h(2^n\mu x)-2\mu h(2^{n-1}x)\|+2^{-n}\|2\mu h(2^{n-1}x)-\mu h(2^nx)\|\\ &\leq &
2^{-n+1}\varphi(2^{n-1}x,2^{n-1}x,0),
\end{eqnarray*}
for all $ \mu\in S^1$ and all $x\in A$.
Since the right side tends to zero as $n\to \infty$, we have $$T(\mu x)=
\displaystyle{\lim_{n\to\infty}}\frac{h(2^n\mu x)}{2^n}=\displaystyle{\lim_{n\to\infty}}\frac{\mu h(2^nx)}{2^n}=\mu T(x),$$
for  all $ \mu\in S^1$ and all  $ x\in A$.
Obviously, $T(0x)=0=0T(x)$.

Next, let $\lambda\in {\bf C} (\lambda\neq 0)$ and let $M$ be a natural number greater than $4|\lambda|$.
Then $|\frac{\lambda}{M}|<\frac{1}{4}<1-\frac{2}{3}=1/3$. By Theorem 1 of [8], there exist three numbers $\mu_1,
\mu_2, \mu_3\in S^1$ such that $3\frac{\lambda}{M}=\mu_1+\mu_2+\mu_3$. By the additivity of $T$
we get $T(\frac{1}{3}x)=\frac{1}{3}T(x)$ for all $ x\in A$. Therefore,
\begin{eqnarray*}
T(\lambda x)& = & T(\frac{M}{3}\cdot 3 \cdot \frac{\lambda}{M}x)=MT(\frac{1}{3}\cdot 3\cdot \frac{\lambda}{M}x)
=\frac{M}{3}T(3\cdot \frac{\lambda}{M}x)\\ & = & \frac{M}{3}T(\mu_1x+\mu_2x+\mu_3x)
=\frac{M}{3}(T(\mu_1x)+T(\mu_2x)+T(\mu_3x)) \\ & = & \frac{M}{3}(\mu_1+\mu_2+\mu_3)T(x)
=\frac{M}{3}\cdot 3\cdot \frac{\lambda}{M}=\lambda T(x),
\end{eqnarray*}
for all $x\in A$.
So that $T$ is ${\bf C}$-linear.

Put $x=y=0$ and replace $z$ by $2^nz$ in {\rm (2)}. Then
$$\frac{1}{2^{3n}}\|h(2^{3n}zz^*z)-h(2^nz)h(2^nz)^*h(2^nz)\|\leq\frac{1}{2^{3n}}\varphi(0, 0, 2^nz)\leq\frac{1}{2^n}
\varphi(0, 0, 2^nz),$$ for all $z\in A$.
Hence \begin{eqnarray*} T(zz^*z) &=& \displaystyle{\lim_{n\to\infty}}\frac{h(2^{3n}zz^*z)}{2^{3n}}
=\displaystyle{\lim_{n\to\infty}}\frac{h(2^nz)}{2^n}\frac{h(2^nz)^*}{2^n}\frac{h(2^nz)}{2^n} \\ & =& T(z)T(z)^*T(z),
\end{eqnarray*}
for all $z \in A$.
It follows that $T$ is a $J^*$-homomorphism satisfying the required inequality. $\Box$

{\bf Example 2.3.} Let $H : A \to A$ be a $J^*$-homomorphism, and let $h : A \to A$ be defined by
$h(x)=\left \{ \begin{array}{cc}H(x) &\qquad  \|x\|<1\\ 0  & \qquad \|x\|\geq 1 \end{array} \right .$ and $\varphi(x, y, z)=4$.

Then
$$\widetilde{\varphi}(x,y,z)=\frac{1}{2}\displaystyle{\sum_{n=0}^\infty}2^{-n}\cdot 4=4,$$ and
$$\|h(\mu x+\mu y+zz^*z)-\mu h(x)-\mu h(y)-h(z)h(z)^*h(z)\|\leq 4=\varphi(x, y, z),$$
for all $\mu\in S^1$ and  all $x, y, z\in A$.

Since in a $C^*$-algebra $$\|xx^*x\|=\|xx^*xx^*xx^*\|^{\frac{1}{2}}=\|(xx^*)^3\|^{\frac{1}{2}}
=\|xx^*\|^{\frac{3}{2}}=\|x\|^{2\cdot \frac{3}{2}}=\|x\|^3,$$
 and every $J^*$-homomorphism is norm decreasing [6], we conclude that $h(xx^*x)=h(x)h(x^*)h(x)$.
Note also that $h$ is not linear. Further, $T(0)=\displaystyle{\lim_{n\to\infty}}\frac{h(0)}{2^n}=0$ and
for $x\neq 0$ we have $$T(x)=\displaystyle{\lim_{n\to\infty}}\frac{h(2^nx)}{2^n}
=\displaystyle{\lim_{n\to\infty}}\frac{0}{2^n}=0,,$$ since for sufficiently large $n, \|2^nx\|\geq 1$. Thus $T$ is identically zero and
 $$\|h(x)-T(x)\|\leq \widetilde{\varphi}(x, x, 0)=4,$$ for all $ x\in A$.

{\bf Corollary 2.4.} {\it Suppose that $h : A \to B$ is a mapping with $h(0)=0$ for which there exist constants $\alpha \geq 0$ and $p \in[0, 1)$ such that
$$
\|h(\mu x+\mu y+zz^*z)-\mu h(x)-\mu h(y)-h(z)h(z)^*h(z)\| \leq  \alpha(\|x\|^p+\|y\|^p+\|z\|^p)
,$$
for all $\mu\in S^1$ and all $ x, y, z\in A.$
Then there is a unique $J^*$-homomorphism $T : A \to B$ such that $$\|h(x)-T(x)\|\leq\frac{\alpha}{1-2^{p-1}}\|x\|^p,$$
for all $ x\in A$.}

{\bf Proof.} Put $\varphi(x, y, z)=\alpha(\|x\|^p+\|y\|^p+\|z\|^p)$ in Theorem 2.2. $\Box$

{\bf Theorem 2.5.} {\it Suppose that $h : A \to B$ is a mapping with $h(0)=0$ for which there exists a function
$\varphi: A\times A\times A\to [0, \infty)$ such that
$$\widetilde{\varphi}(x, y, z): =\frac{1}{2}\displaystyle{\sum_{n=0}^\infty}2^{-n}\varphi(2^nx,2^ny,2^nz)<\infty,$$
\begin{equation}
\|h(\mu x+\mu y+zz^*z)-\mu h(x)-\mu h(y)-h(z)h(z)^*h(z)\|\leq \varphi(x,y,z),
\end{equation}
for $\mu=1,i$ and all $x, y, z\in A$. If for each fixed $x\in A$ the function $t\mapsto h(tx)$ is continuous
 on ${\bf R}$, then there exists a unique $J^*$-homomorphism $T : A \to B$ such that
$$\|h(x)-T(x)\|\leq\widetilde{\varphi}(x, x, 0),$$ for all $ x\in A$.}

{\bf Proof.} Put $z=0$ and $\mu=1$ in {\rm (3)}. It follows from the G\u avruta theorem that
there exists a unique additive mapping $T : A \to B$ given by $$T(x)=\displaystyle{\lim_{n\to\infty}}\frac{h(2^nx)}{2^n},$$
for all $ x\in A$.

By the same reasoning as in the proof of the main theorem of [11], the mapping $T$ is ${\bf R}$-linear.

Assume $y=z=0$ and $\mu=i$. It follows from {\rm (3)} that $$\|h(ix)-ih(x)\|\leq\varphi(x, 0, 0),$$
 for all $ x\in A$. Hence $$\frac{1}{2^n}\|h(2^nix)-ih(2^nx)\|\leq\varphi(2^nx, 0, 0) ,$$ for all $x \in A$.
 The right side tends to zero as $n\to\infty$, so  $$T(ix)=\displaystyle{\lim_{n\to\infty}}\frac{h(2^nix)}{2^n}
=\displaystyle{\lim_{n\to\infty}}\frac{ih(2^nx)}{2^n}=iT(x),$$ for all $ x\in A$.
For every $\lambda\in {\bf C}, \lambda=s+it$ in which $s,t\in{\bf R}$ we have
\begin{eqnarray*}
T(\lambda x) & = & T(sx+itx)=sT(x)+tT(ix)=sT(x)+itT(x)=(s+it)T(x)
\\ & = & \lambda T(x),
\end{eqnarray*} for all $x \in A$.
Thus $T$ is ${\bf C}$-linear. $\Box$

Choonkil Baak\\
Department of Mathematics, Chungnam National University, Daejeon 305-764, South Korea\\
e-mail: cgpark@cnu.ac.kr \\
Mohammad Sal Moslehian\\
Department of Mathematics, Ferdowsi University, P. O. Box 1159, Mashhad 91775, Iran\\
e-mail: moslehian@ferdowsi.um.ac.ir
\end{document}